\documentclass{amsart}
\usepackage{hyperref}
\usepackage{algpseudocode}

\hypersetup{
    pdftitle={Algorithm xxx: A Correctly Rounded Newton Step for the Reciprocal Square Root},
    pdfauthor={Carlos F. Borges},
    pdfkeywords={rsqrt(), floating-point, fused multiply-add, IEEE 754}
    }

\newtheorem{algorithm}{Algorithm}
\numberwithin{equation}{section}

\def\F{\mathbb{F}}
\def\Re{\mathbb{R}}
\DeclareMathOperator{\ulp}{ulp}
\DeclareMathOperator{\bu}{\mathbf{u}}

\begin{document}

\title{A Correctly Rounded Newton Step for the Reciprocal Square Root}

\author{Carlos F. Borges}
\address{\\Department of Applied Mathematics\\Naval Postgraduate School\\Monterey CA 93943}
\email{borges@nps.edu}
        
\begin{abstract}
The reciprocal square root is an important computation for which many sophisticated algorithms exist (see for example \cite{Moroz,863046,863031} and the references therein). A common theme is the use of Newton's method to refine the estimates. In this paper we develop a correctly rounded Newton step that can be used to improve the accuracy of a naive calculation (using methods similar to those developed in \cite{borges}) . The approach relies on the use of the fused multiply-add (FMA) which is widely available in hardware on a variety of modern computer architectures. We then introduce the notion of {\em weak rounding} and prove that our proposed algorithm meets this standard.  We then show how to leverage the exact Newton step to get a Halley's method compensation which requires one additional FMA and one additional multiplication. This method appears to give correctly rounded results experimentally and we show that it can be combined with a square root free method for estimating the reciprocal square root to get a method that is both very fast (in computing environments with a slow square root) and, experimentally, highly accurate.
\end{abstract}




\date{April 30, 2019.}



\maketitle

There are many current algorithms (see for example \cite{Moroz,863046,863031} and the references therein) in the literature for computing the reciprocal square root. In this paper we show how one can leverage the fused multiply-add to compute a correctly rounded Newton step, that will yield a very accurate answer. We then show how to leverage the correctly rounded Newton step to get a Halley's method compensation which requires one additional FMA and one additional multiplication. When this approach is combined with a square root free method for estimating the reciprocal square root we get a method that is both very fast and, experimentally, highly accurate. Such methods are important in computing environments that do not have a fast square root (e.g. microcontrollers and FPGAs). We illustrate the accuracy of the algorithms {\em experimentally} by comparing them to extended precision calculations carried out with the MPFR package in Julia. All of these algorithms can be implemented on architectures that do not have a hardware FMA by using a software implementation of the FMA although this would likely signficantly reduce the speed.

This paper will be restricted to the case where all floating-point calculations are done in IEEE 754 compliant radix~2 arithmetic using round-to-nearest rounding, although many of the results can be extended to other formats under proper conditions.

\section{Mathematical Preliminaries}

We begin with a few definitions. We shall denote by $\F \subset \Re$ the set of all radix 2 floating-point numbers with precision $p$. We define $fl(x):\Re \rightarrow \F$ to be a function such that $fl(x)$ is the element of $\F$ that is closest to $x$. Since we have restricted ourselves to radix 2 IEEE 754 compliant arithmetic we will assume that round-to-even is used in the event of a tie. Throughout this paper we assume that $ulp(x)$ is a unit in the last place as defined in \cite{Goldberg}. To wit, we define $ulp(x) = 2^{e-p+1}$ for any $x \in [2^e,2^{e+1})$ where $p$ is the precision of the floating-point format. For example, in double precision $\ulp(1) = 2^{-52}$. We define the {\em unit roundoff}, which we will denote by $\bu = \frac{1}{2}\ulp(1)$ so that, for example, in double precision $\bu = 2^{-53}$.

\section{Computing the reciprocal of the square root}

Given a floating-point number $x \in \F$ we wish to compute $$y = \frac{1}{\sqrt{x}}$$ in floating-point. The most accurate naive approach is simply\footnote{One can also use y=1/sqrt(x) but this can lead to errors greater than 1 ulp. See \cite{Markstein2}.}

\begin{algorithm} Naive rsqrt - {\tt rsqrtNaive(x)}

\hrule
\begin{algorithmic}
	\State {r = 1/x}
	\State {y = sqrt(r)}
\end{algorithmic}
\hrule
\label{Naiversqrt}
\end{algorithm}

We note that the computed quantity $y$ is subject to various errors due to the effects of finite precision and hence
\begin{equation}
\sqrt{\frac{1}{x}} = y (1+\nu)
\label{recip1}
\end{equation}
for some $\nu \in \Re$, where it is reasonable to assume that $|\nu| < 1$. Under our assumptions as to the floating-point environment, a standard error analyis reveals that 
\begin{equation}
|\nu| < \frac{3}{2}\bu + O(\bu^2).
\label{crap}
\end{equation}

Squaring both sides of \ref{recip1} and a bit of algebra gives us
\begin{equation}
x y^2 = \frac{1}{(1+\nu)^2}
\label{recip2}
\end{equation}
and then two steps of long division on the right hand side and a bit more algebra yields
\begin{equation}
\nu = \frac{1-x y^2}{2}  + \nu^2\frac{3+2\nu}{2(1+\nu)^2}.
\label{recip3}
\end{equation}

Observe that at this point that for any $y$ satisfying \ref{recip1} with $|\nu| < 1$ we have
$$
\sqrt{\frac{1}{x}}  =  y + y\left(\frac{1-xy^2}{2} + \nu^2\frac{3+2\nu}{2(1+\nu)^2}\right) 
$$
and therefore
\begin{equation}
y + y\frac{1-xy^2}{2} = \sqrt{\frac{1}{x}} + O(y\nu^2).
\label{bigone}
\end{equation}
If we ignore the $O(y\nu^2)$ term and let
\begin{equation}
\bar{\nu} = \frac{1-xy^2}{2}.
\label{nu}
\end{equation}
We can add a compensation to the naively computed value using
\begin{equation}
y_C = y + y\bar{\nu}
\label{Newtonstep}
\end{equation}
which the astute reader will recognize as the Newton iteration for
$$
f(y) = x - \frac{1}{y^2}.
$$
This is a common approach to estimating the reciprocal square root (see \cite{Moroz} and references therein). We note that the point of this derivation is not to reinvent the Newton iteration but rather to provide the relationship in \ref{bigone} that we will require in our error analysis.

We can perform a traditional rounding error analysis for the single Newton step in \ref{Newtonstep} using the standard model for radix 2 round-to-nearest floating point arithmetic. Assume that the computed value of $\bar{\nu}$ is  $\bar{\nu}(1+\epsilon)$ for some $|\epsilon| < 1$. Then, if we apply the Newton step using an FMA, the final compensated value, $y_C$, involves only a single rounding operation and therefore satisfies
\begin{eqnarray*}
y_C & = &  (y + y  \frac{1-xy^2}{2}(1+\epsilon))(1+\delta_1) \\
& = &  (y + y  \frac{1-xy^2}{2})(1+\delta_1) + y  \frac{1-xy^2}{2}\epsilon(1+\delta_1)
\end{eqnarray*}
where $|\delta_1| \leq \bu$.Using \ref{bigone} to replace the first term in parentheses on the right we get
$$
y_C =  \sqrt{\frac{1}{x}}(1+\delta_1) + O(y\nu^2) + y \frac{1-xy^2}{2}(\epsilon (1+ \delta_1))
$$
and then using \ref{recip3} gives
$$
y_C =  \sqrt{\frac{1}{x}}(1+\delta_1) + O(y\nu^2) + O(y\nu\epsilon)
$$
And finally
$$
\frac{\left| y_C - \sqrt{\frac{1}{x}}\right|}{\sqrt{\frac{1}{x}}}  = \frac{\left|y\left( \delta_1 + O(\nu^2) + O(\nu \epsilon) \right)\right|}{y(1+\nu)} = \frac{\left|\left( \delta_1 + O(\nu^2) + O(\nu \epsilon) \right)\right|}{(1+\nu)}.
$$
Clearly, if both $\nu$ and $\epsilon$ are $O(\bu)$ then the relative error of our approximation is bounded by
\begin{equation}
\frac{\left| y_C - \sqrt{\frac{1}{x}}\right|}{\sqrt{\frac{1}{x}}}
\leq \bu + O(\bu^2)
\label{goodbound}
\end{equation}
This is a useful condition as it implies that $y_C$ is {\em faithful} and, even more, that the true value can't be much further than halfway across the gap between the computed result and the next (or previous) floating-point number. It's not quite correctly rounded but it's very close.

Now the important part, when $\nu$ is very small computing $\bar{\nu}$ as described in equation \ref{nu} can involve extreme cancellation. As a result, the relative error, $\epsilon$, can be very large and this can spoil our attempt at compensation. To avoid this consider rewriting equation \ref{nu} in the following suggestive form
$$
\bar{\nu} = \frac{1}{2}(1-xr - x(y^2-r))
$$
and using the following Newton's Method compensation which requires adding four FMA calls and a single multiply to the naive algorithm:

\begin{algorithm} Newton's Method Compensation - {\tt rsqrtNewton(x)}
\hrule
\begin{algorithmic}
   \State {r = 1/x}
	\State {y = sqrt(r)}
	\State {mxhalf = -0.5*x}
   \State {$\sigma$ = fma(mxhalf, r, 0.5)}
   \State {$\tau$ = fma(y, y, -r)}
   \State {$\bar{\nu}$ = fma(mxhalf, $\tau$, $\sigma$)}
   \State {y = fma(y, $\bar{\nu}$, y)}
\end{algorithmic}
\hrule
\label{Correctedrsqrt}
\end{algorithm}

\subsection{Error Analysis}

To see why algorithm \ref{Correctedrsqrt} works we need only observe that under appropriate conditions (see Chapter 4 Theorems 4.9 and 4.10 in \cite{FPHB}), both $\sigma, \tau \in \F$,  and are computed exactly using the FMA (see \cite{Markstein} or Chapter 4 Corollary 4.11 and 4.12 in \cite{FPHB} for the precise conditions). This means that the computed value of $\bar{\nu}$ is correctly rounded (using an FMA) and hence $\epsilon = \bu$ which implies that the result of algorithm \ref{Correctedrsqrt} satisfies the bound in equation \ref{goodbound} and hence the algorithm gives a weakly rounded result.

\subsection{Numerical Testing}

Since our error bound is not sufficiently tight to show correct rounding, we now test the compensated algorithm against the naive approach. All testing is done in IEEE 754 double precision arithmetic with code written in Julia 1.5.3 running on an Intel(R) Core(TM) i7-7700K CPU @ 4.20GHz. As a baseline for testing purposes we will compute the reciprocal square root, $\bar{y}$, by using the BigFloat format in Julia which uses the GNU MPFR package to do an extended precision calculation. It is critical to note that neither square root nor division are finite operations. That is, both can result in infinite length results and hence this extended precision computation is necessarily prone to double rounding. That means that the result of the MPFR computation cannot be guaranteed to represent a correctly rounded value of the true result (although it will do so in nearly every case).

 We will do this using $10^{9}$ uniformly distributed double precision (Float64) random inputs. We will run both algorithms on each random input as well as computing $\bar{y}$. In the tables below we summarize the error rates of each algorithm which we define to be the percentage of times each algorithm differed from $\bar{y}$ by exactly zero ulp and exactly one ulp. We note that neither algorithm ever differed from $\bar{y}$ by more than one ulp.
 
\begin{table}[!ht]
  \begin{center}
    \caption{Error rate for computing the reciprocal square root with the Newton Compensation versus the Naive computation.}
    \label{tab:tablersqrt1}
    \begin{tabular}{l l l}
      $x \sim \mathcal{U}(1,1/2)$ & Newton & Naive \\
      \hline
Zero ulp &  100 & 89.227\\
One ulp  & 0 & 10.773 \\
\\
 $x \sim \mathcal{U}(1,2)$& Newton & Naive \\
      \hline
Zero ulp & 100 & 84.762  \\
One ulp  & 0 & 15.238
    \end{tabular}
  \end{center}
\end{table}

Although the compensated algorithm is clearly very accurate, it does not always yield a correctly rounded result. To wit, let $x = 1-2\bu$, then the computed value of $y=1.0$ and the compensation will be computed exactly as $y\bar{\nu} = \bu$, and $1+\bu \rightarrow 1$ because of round-to-even.\footnote{We are not aware of any other examples beyond this one and any multiples of the form $(1-2\bu)4^k$.} Since we know that the reciprocal square root of any radix 2 floating point number cannot be the exact midpoint of two consecutive floating point numbers (see Chapter 4 Theorem 4.20 in \cite{FPHB}) it is clear that the problem in this case is the difference between $\bar{\nu}$ and $\nu$.

Note that  $\bar{\nu}$ and $\nu$ always have the same sign, and further that
\begin{equation}
\bar{\nu} < \nu
\label{overcorrect}
\end{equation}
provided that $\nu \neq 0$,which is seen by rearranging \ref{recip3}. This means that algorithm \ref{Correctedrsqrt} will always slightly undercompensate if $\bar{\nu}$ is positive, and slightly overcompensate if it is negative. Now, in the specific example just discussed, the compensation was positive and hence too small and we really should have compensated by a tad more which would have given us the correctly rounded result of $y=1+2\bu$. We note that although the result gend{equation}rated by the algorithm is not correctly rounded, it is {\em weakly rounded} as guaranteed by the earlier analysis.

\subsection{Halley's Method Compensation}

The fact that we can compute a correctly rounded $\bar{\nu}$ means it may reasonably be used to estimate $\nu$ directly. To that end, if $\nu = O(\bu)$ then equation \ref{recip3} can be rewritten as
$$
\nu = \bar{\nu}  + \frac{3}{2}\nu^2 + O(\bu^3).
$$
A little algebra introduces some more $O(\bu^3)$ terms and yields
\begin{eqnarray*}
\nu  & = & \frac{\bar{\nu}}{1 - \frac{3}{2}\nu} + O(\bu^3)\\
& = & \bar{\nu} \left(1 + \frac{3}{2}\nu \right) + O(\bu^3)
\end{eqnarray*}
and finally, replacing $\nu$ on the right hand side with $\bar{\nu}$ introduces yet more $O(\bu^3)$ terms and gives
\begin{equation}
\nu  = \bar{\nu} \left(1 + \frac{3}{2}\bar{\nu} \right) + O(\bu^3).
\label{ohyeah}
\end{equation}
The careful reader will notice that this is simply Halley's method.

\begin{algorithm} Halley's Method Compensation - {\tt rsqrtHalley(x)}
\hrule
\begin{algorithmic}
   \State {r = 1/x}
	\State {y = sqrt(r)}
	\State {mxhalf = -0.5*x}
   \State {$\sigma$ = fma(mxhalf, r, 0.5)}
   \State {$\tau$ = fma(y, y, -r)}
   \State {$\bar{\nu}$ = fma(mxhalf, $\tau$, $\sigma$)}
    \State {$\nu$ = fma(1.5*$\bar{\nu}$,$\bar{\nu}$,$\bar{\nu}$)}
   \State {y = fma(y, $\nu$, y)}
\end{algorithmic}
\hrule
\label{Correctedrsqrt2}
\end{algorithm}

This algorithm requires one more multiply and one more FMA than the Newton's method compensation. In extensive testing this algorithm has always gend{equation}rated the correctly rounded result, even for $x = 1-2\bu$. The author conjectures that this algorithm always gives a correctly rounded result but does not have a proof.

\subsection{Developing a fast square-root free variant}

The compensation can be applied to existing square-root-free algorithms for computing the reciprocal square root. Such methods are important in computing environments that do not have a fast square root (e.g. microcontrollers and FPGAs). There are many such algorithms but we have chosen to use one of a type that is best known for appearing in the code for the video game {\em Quake III Arena}. The specific algorithm from \cite{Moroz} that we will investigate is called {\tt RcpSqrt331d(x)}. A Julia port of that code appears below.

\begin{verbatim}
function RcpSqrt331d(x::Float64)
    i = reinterpret(Int64,x)
    k = i & 0x0010000000000000
    if k != 0
        i = 0x5fdb3d14170034b6 - (i >> 1)
        y = reinterpret(Float64,i)
        y = 2.33124735553421569*y*fma(-x, y*y,1.07497362654295614 )
    else
        i = 0x5fe33d18a2b9ef5f - (i >> 1)
        y = reinterpret(Float64,i)
        y = 0.82421942523718461*y*fma(-x, y*y, 2.1499494964450325)
    end
    mxhalf = -0.5*x
    y = y*fma(mxhalf, y*y, 1.5000000034937999)
  # The next two lines are a single step of Newton
    r = fma(mxhalf, y*y, 0.5)
    y = fma(y, r, y)
end
\end{verbatim}

Note that the final two lines of the code are simply a careful application of one step of Newton which requires two FMAs and one multiply. We replace those lines with the correctly rounded Newton's method compensation. The added lines are:

\begin{algorithm} {\tt RcpSqrt331dHalley(x)} - partial code
\hrule
\begin{algorithmic}
	\State{\# The following lines replace the Newton step in RcpSqrt331d(x)}
   \State {r = 1/x}
   \State {$\sigma$ = fma(r,mxhalf,.5)}
   \State {$\tau$ = fma(y, y, -r)}
   \State {$\bar{\nu}$ = fma(mxhalf, $\tau$, $\sigma$)}
    \State {$\nu$ = fma(1.5*$\bar{\nu}$,$\bar{\nu}$,$\bar{\nu}$)}
   \State {y = fma(y, $\nu$, y)}
\end{algorithmic}
\hrule
\label{FastCorrectedrsqrt}
\end{algorithm}

This adds three FMAs and one divide to the code and benchmark timings in Julia indicate that it adds about 5\% to the execution time.

In the table below we show the results from the accuracy test comparing the two algorithms. Note that the compensated form returns the correctly rounded answer every single time in the experiment. Moreover, it gives a correctly rounded answer for $x = 1-2\bu$.

\begin{table}[!ht]
  \begin{center}
    \caption{Error rate for computing the reciprocal square root with Halley compensation in the square root free algorithm versus the uncompensated square root free algorithm.}
    \label{tab:tablefrsqrt1}
    \begin{tabular}{l l l}
      $x \sim \mathcal{U}(1/2,1)$ & RcpSqrt331dHalley & RcpSqrt331d \\
      \hline
      Zero ulp & 100 & 87.324\\
      One ulp & 0 & 12.676 \\
      \\
      $x \sim \mathcal{U}(1,2)$ & RcpSqrt331dHalley & RcpSqrt331d \\
      \hline
      Zero ulp &  100 & 82.119\\
      One ulp  & 0 & 17.881 
    \end{tabular} 
  \end{center}
\end{table}

\section{Conclusions}

We have shown that it is possible to compute the correctly rounded value of the standard Newton step for reciprocal square root for estimates $y$ that are sufficiently close to the true value. This Newton step can be used to get a compensated estimate of the reciprocal square root that is weakly rounded which is a strong guarantee of high accuracy. We then showed how to extend this to a Halley's method compensation. This leads to an algorithm that we conjecture is correctly rounded although we cannot prove it at this time. Moreover, these compensation methods can be shown to significantly improve the accuracy of certain square-root free methods for estimating the reciprocal square root. 

\bibliographystyle{acm}

\end{document}